\documentclass[reqno]{amsart}
\usepackage{amsfonts}
\usepackage{amssymb}
\usepackage{amsthm}
\usepackage{upref}
\usepackage{enumerate}
\makeatletter
\def\LaTeX{\leavevmode L\raise.42ex
    \hbox{\kern-.3em\size{\sf@size}{0pt}\selectfont A}\kern-.15em\TeX}
 \makeatother
 
\numberwithin{equation}{section}

\newtheorem{lemma}{Lemma}[section]
\newtheorem{theorem}[lemma]{Theorem} 
\newtheorem{corollary}[lemma]{Corollary}
\newtheorem{proposition}[lemma]{Proposition}

\theoremstyle{definition}

\newtheorem{assumption}[lemma]{Assumption}

\newtheorem{remark}[lemma]{Remark}

  \newcommand{\e}{\eqref}
\newcommand{\ri}{\rightarrow}
\newcommand{\q}{\quad}
\newcommand{\ii}{\infty}
\newcommand{\h}{\hbar}
\renewcommand{\d}{\delta}

    \newcommand{\A}{\operatorname{Ai\,}}
     \newcommand{\B}{\operatorname{Bi\,}}

\newenvironment{pf}{\begin{proof}}{\end{proof}}

\def\qqq{\mathrel{\subset\mkern-15mu\lower.38ex\hbox{${\scriptscriptstyle\rightarrow}$}}}

\let\Bbb\mathbb

\begin{document}
\title[semiclassical limit of eigenfunctions]
{The semiclassical limit of eigenfunctions of the Schr\"odinger equation and the Bohr-Sommerfeld quantization condition, revisited}
\author{ D. R. Yafaev  }
\address{ IRMAR, Universit\'{e} de Rennes I\\ Campus de
  Beaulieu, 35042 Rennes Cedex, FRANCE}
\email{yafaev@univ-rennes1.fr}
\thanks{Partially supported by the project NONAa, ANR-08-BLANC-0228}
 
\dedicatory{To Vasilij Mikhailovich Babich on his 80-th birthday}

\begin{abstract}
Consider the semiclassical  limit, as the Planck constant $\hbar\ri 0$, of     bound  states of a   quantum particle in a one-dimensional potential well. We justify
   the semiclassical  asymptotics  of   eigenfunctions and recover the Bohr-Sommerfeld quantization condition. 
   \end{abstract}
\maketitle

\thispagestyle{empty}

\section{Introduction}
{\bf 1.1. }
We study  the limit as $\h\ri 0$ of   eigenfunctions $\psi (x )=\psi (x ;\lambda,\h) $ of the Schr\"odinger equation
\begin{equation}
 -\h^2  \psi ''(x ) +v(x) \psi (x ) =\lambda \psi (x ), \q v(x)=\overline{v(x)}, \q \psi \in L^{2}({\Bbb R} ),
\label{eq:I2}\end{equation}
for   $\lambda $   close to some non-critical energy $\lambda_{0}$ (that is $  v'(x) \neq 0$ for   $x$ such that $v(x)=\lambda_{0}$). We   assume that the equation
  $v(x)=\lambda $ has exactly two solutions (the turning points) $x_{\pm}=x_{\pm}(\lambda )$ and that $v(x)<\lambda $ for $x\in (x_{-}, x_{+})$. Thus,  $(x_{-}, x_{+})$ is a potential well and the energy $\lambda $ is separated from its bottom.
We   suppose that eigenfunctions $\psi (x )$ are real and normalized, that is
 \[
  \int _{-\infty}^\infty   \psi ^2(x )  dx =1.
  \]

It is a common wisdom that the limit of  $\psi (x )=\psi (x ;\lambda,\h) $ as $\h\ri 0$ is described by the Green-Liouville approximation away from the turning points $x_{\pm} $. In neighborhoods of the turning points the asymptotics of $ \psi (x )$ is more complicated and is given in terms of an Airy function. Surprisingly, we have not found  a precise formulation and a proof of this result in the literature. Our goal is to fill in this gap. We follow here the scheme suggested by R.~E.~Langer and thoroughly exposed   by F.~W.~Olver in his book \cite{Olver}.

The detailed asymptotics of $ \psi (x )$ described in Theorems~\ref{Airy}  and \ref{norm} allows one to recover the classical Bohr-Sommer\-feld quantization condition on $\lambda $ (see Theorem~\ref{pm}). Actually, we prove somewhat more establishing a one-to-one correspondence between eigenvalues of the Schr\"odinger operator 
$H_{\h}= -\h^2  d^2/dx^2+v(x) $ from a neighborhood of a non-critical energy and points $(n+1/2)\h$ where $n$ is an integer. This implies the semiclassical Weyl formula for the distribution of eigenvalues of the   operator 
$H_{\h}  $ as $\h\to 0$  with a strong estimate of the remainder. It turns out (see Corollary~\ref{pmc}) that  this remainder never exceeds $1$.  We also obtain in Theorem~\ref{disc} the quantization condition for discontinuous functions $v(x)$. This formula generalizes that of Bohr and Sommerfeld and is probably new.

We note that the Bohr-Sommer\-feld quantization condition is  also  well known in much more difficult multidimensional problems. In this context we mention book \cite{BB}  by V.~M.~Ba\-bich and V.~S.~Buldyrev (where the ray approximation  is used), book   \cite{FedMas}  by M.~V.~Fedoryuk and V.~P.~Maslov (where the Maslov canonical operator is used) as well as papers \cite{ShRo} by  B. Helffer  et D. Robert and
   \cite{ShMaRo} by B. Helffer, A. Martinez  and D.~Robert (where the methods of microlocal analysis are used).
   
    However, in the one-dimensional problem it is more natural to rely on methods of ordinary differential equations. Such an approach was developed 
  by M.~V.~Fedoryuk (see his book \cite{Fed}) for analytic potentials. In this case one can avoid a study of turning points so that the Airy function does not appear.

 \medskip
 
 {\bf 1.2. }
   The asymptotics of   eigenfunctions     yields (see Proposition~\ref{normw}) asymptotics of  observables
   \begin{equation}
 \int _{-\infty}^\infty w(x)\psi ^2(x; \lambda,\h)dx 
\label{eq:Obs}\end{equation}
for sufficiently arbitrary functions $w(x)$.  For example, we can take  for $w(x) $ characteristic functions of Borel subsets of ${\Bbb R} $ or choose $w(x)=v(x)$. This gives the asymptotics of the kinetic energy
 \begin{equation}
K(\lambda, \h):= \h^2   \int _{-\infty}^\infty     \psi'(x ; \lambda,\h )^2 dx = K_{cl}(\lambda ) + O(\h^{1/3 }) 
 \label{eq:mainkb}\end{equation}
  as $\h\ri 0$ uniformly
for $\lambda$ from a neighborhood of the point $\lambda_{0}$.  The leading term $K_{cl}(\lambda )$ (the index $``cl$" stands of course  for the corresponding classical object) is given by the expression
\begin{equation}
K_{cl}(\lambda  )=
  \int _{x_{-}(\lambda )}^{x_{+}(\lambda )}   ( \lambda  -v(x))^{1/2} dx 
\Bigl(   \int _{x_{-}(\lambda )}^{x_{+}(\lambda )}(\lambda  -v(x))^{-1/2} dx \Bigr)^{-1}.
\label{eq:mainOld}\end{equation}
Note that the integrals here are taken over 
  the classically allowed region and  that $K_{cl}(\lambda )$ coincides (see subsection~4.3) with the averaged value  of the kinetic energy of a particle of energy $\lambda$ in classical mechanics.

 We emphasize that our derivation of the quantization condition and of asymptotic formulas for observables \e{eq:Obs} requires Airy functions although they do not enter into the final answer. However,  we do not know  how to avoid Airy functions without additional assumptions on $v(x)$.

\section{Semiclassical solutions of the Schr\"odinger equation}


 \medskip
 
 {\bf 2.1.}
   It is convenient to rewrite equation \e{eq:I2} as
\begin{equation}
 - u_{\h}^{\prime \prime}(x ) +\h^{-2}  q(x) u_{\h}(x ) =0,
\label{eq:R1}\end{equation}
where
\[
q(x)= q(x;\lambda)= v(x)-\lambda.
\]

We need some regularity of the function $v(x)$ and a weak condition on its behavior   at infinity.

 \begin{assumption}\label{infi}
 The function $v\in C^2 ({\Bbb R})$  and, for some $\rho_{0}> 1$, the function
 \[
 \big(| q (x)|^{-3}  q^{\prime  } (x)^2 + q (x)^{-2} | q^{\prime \prime} (x)| \big)
\Big|\int_0^x |q(y)|^{1/2} dy \Big|^{\rho_{0}}
\]
 is bounded  for  sufficiently large $|x|$.
 \end{assumption}

 The last condition is satisfied   in all reasonable cases. For example, if $v(x)\ri v_{0}> \lambda$,  it is sufficient to require that
  \[
  v^{\prime  } (x) ^2+   | v^{\prime \prime} (x)|=O( |x|^{-\rho_{0}}),\q \rho_{0}> 1 ,\q |x|\to \infty.
\]
It is also satisfied if $v(x)$ behaves at infinity as $|x|^\alpha$ or $e^{\alpha|x|}$ where $\alpha>0$; in these cases $\rho_{0}=2$.

We consider   the case of one potential well. To be more precise, we make the following

 \begin{assumption}\label{onewell}
 The equation $v(x)=\lambda$ has two solutions $x_{+} = x_{+} (\lambda)$ and $x_{-}= x_{-} (\lambda)$. We suppose that    $x_{-} < x_{+}$, $v  (x ) < \lambda$ for $x\in (x_{-}, x_{+})$,  $v  (x ) > \lambda$  for $x\not\in [x_{-}, x_{+}]$  and
  \[
\liminf_{|x|\ri\ii}v(x)>\lambda.
\]
Moreover, the function $v$ belongs to the class $C^3$ in some  neighborhoods of the points $x_{\pm}$ and $\pm v' (x_{\pm}) > 0$.
  \end{assumption}

 Note that if Assumption~\ref{onewell} is satisfied for some $\lambda_{0}$, then it is also satisfied for all $\lambda$ from some neighborhood of $\lambda_{0}$.

Our   goal in this section is to describe asymptotics  as $\h\ri 0$ of   solutions  $u_+  (x )= u_+  (x; \lambda,\h )$ and $u_-  (x )= u_-  (x; \lambda,\h )$ of equation
\e{eq:R1} exponentially decaying as $x \ri +\ii$ and $x \ri -\ii$, respectively. These asymptotics will be given in terms of  an  Airy function and are uniform with respect to $x\in [x_{1},\ii)$ or $x\in (-\ii, x_{1}]$ where  $x_{1}$  is an arbitrary point from the interval  $ (x_{-}, x_{+})$.

 \medskip
 
 {\bf 2.2. }
 Let us recall the definition of   Airy functions and their necessary properties (see, e.g., \cite{Olver}), for details). Consider   the equation
\begin{equation}
 - w^{\prime \prime}(t) +t w (t) =0 
  \label{eq:R5}\end{equation}
  and denote by $\A(t)$ its solution with asymptotics
\begin{equation}
 \A(t) = 2^{-1} \pi^{-1/2} t^{-1/4}  \exp (-2t^{3/2}/3) (1+ O(t^{-3/2})), \q t\ri +\ii.
\label{eq:R6}\end{equation}
Then
\begin{equation}
 \A(t) = \pi^{-1/2} |t|^{-1/4}   \sin (2 |t|^{3/2}/3 +  \pi/4) + O(|t|^{-7/4}) , \q t\ri -\ii.
\label{eq:R7}\end{equation}
Note   that $ \A(t)> 0$ for all $t\geq 0$.

The solution $\B(t)$ of equation \e{eq:R5} is defined by its asymptotics as $t\ri -\ii$ which differs from \e{eq:R7} only by the phase shift: 
\begin{equation}
 \B(t) = - \pi^{-1/2} |t|^{-1/4}   \sin (2 |t|^{3/2}/3- \pi/4) + O(|t|^{-7/4}) , \q t\ri -\ii.
\label{eq:R7B}\end{equation}
For $t\geq 0$, this function is positive and satisfies the estimate
\begin{equation}
 \B(t) \leq C(1+ t)^{-1/4}  \exp (2t^{3/2}/3)  .
\label{eq:R6B}\end{equation}
Here and below we denote by $C$ and $c$ different positive constants whose precise values are of no importance.

We   also  use that all asymptotics \e{eq:R6}, \e{eq:R7} and \e{eq:R7B} can be differentiated in  $t$. In particular, 
the Wronskian
\[
 \{ \A (t) , \B_{-}(t)\} : = \A ^\prime (t)  \B(t) - \A (t)   \B^\prime(t)= -\pi ^{-1}.  
\]
It follows that
\begin{equation}
  \B (s) \A ^{-1}(s)   - \B (t) \A ^{-1}(t)  
 =  \pi ^{-1} \int_{t}^s  \A^{-2}(\tau) d \tau, \q s\geq t\geq 0.
\label{eq:R7cc}\end{equation}

 \medskip
 
 {\bf 2.3. }
 To formulate   results, we need the following auxiliary functions $\xi_{\pm}(x)= \xi_{\pm}(x; \lambda)$: 
\begin{equation}
\begin{split}
\xi_{+}(x)&= \Bigl( \frac{3}{2}\int_{x_{+}}^x q (y)^{1/2} dy\Bigr)^{2/3}, \q x\geq x_{+},
 \\
\xi_{+}(x)&= - \Bigl( \frac{3}{2}\int^{x_{+}}_x |q (y)|^{1/2} dy \Bigr)^{2/3}, \q 
x_{-} < x\leq x_{+},
\end{split}
\label{eq:R3}\end{equation}
and
\[
\begin{split}
\xi_{-}(x) &= \Bigl( \frac{3}{2}\int^{x_{-}}_x q (y)^{1/2} dy \Bigr)^{2/3}, \q x\leq x_{-},
 \\
\xi_{-}(x)&= - \Bigl( \frac{3}{2}\int_{x_{-}}^x |q (y)|^{1/2} dy \Bigr)^{2/3}, \q 
x_{-}\leq x < x_{+}.
\end{split}
\]
Here is a list of properties of these functions. The following result is practically the same as Lemma~3.1 from Chapter~11 of \cite{Olver}.

 \begin{lemma}\label{xi}
  Let $x_{1}\in (x_{-},x_{+})$. Then  $ \xi_{+}\in C^3 (x_{1}, \infty)$,   $ \xi_{-}\in C^3(-\infty, x_{1} )$ and $\xi_{\pm} (x) \to + \ii$ as $x\ri\pm\ii$. The derivatives  
  \begin{equation}
  \pm \xi_{\pm}'  (x) >0, \q  \xi_{\pm}'  (x_{\pm}) =\pm |v' (x_{\pm})|^{1/3}
  \label{eq:R4x}\end{equation}
   and the functions  $\xi_{\pm}  (x) $ satisfy the equation
\begin{equation}
\xi_{\pm}  ^\prime(x)^2 \xi_{\pm}  (x)= q(x).
\label{eq:Req}\end{equation} 
 \end{lemma}
 
 It follows from this lemma that the function 
 \begin{equation}
   p_{\pm} ( x )= ( |\xi_{\pm}^\prime(x)|^{-1/2})^{\prime \prime}| \xi_{\pm}^\prime(x)|^{-3/2} 
\label{eq:Ra2}\end{equation}
 is continuous. Moreover, using   identity \e{eq:Req}, we see that
 \begin{equation}
 -16 p_{\pm} ( x )= 5 \xi_{\pm} (x)^{-2}+ \xi_{\pm} (x)  \big(4 q (x)^{-2}  q^{\prime\prime}  (x)-5 q (x)^{-3} q^{\prime}  (x)^2 \big), \q x \neq x_{\pm},   
\label{eq:Pp}\end{equation}
and hence according to Assumption~\ref{infi}  
 \begin{equation}
  |p_{\pm} ( x )|\leq C | \xi_{\pm} ( x)|^{- 1/2-\rho  },\q  \q  \rho= 3 \min\{\rho_{0}  - 1, 1\} /2> 0.
\label{eq:Ra2x}\end{equation}

 \medskip
 
 {\bf 2.4. }
Let us construct solutions  $u_\pm  (x )= u_\pm (x; \lambda,\h )$   of equation \e{eq:R1} with semiclassical asymptotics as $\h\to 0$ or (and) $x\to \pm\infty$. We define these solutions by their asymptotics as $x\to \pm\infty$. Below all asymptotic relations are supposed to  be differentiable with respect to $x$. In this subsection, we only formulate results.

 \begin{proposition}\label{inf}
Under Assumption~\ref{infi} for every fixed $\h>0$, equation  \e{eq:R1} has a $($unique$)$ solution $u_\pm (x )$ such that
   \begin{align*}
u_\pm  (x )= &2^{-1}  \pi^{1/2} \h^{1/6} q(x)^{-1/4} \exp\Big(\mp \h^{-1} \int_{x_{\pm}}^x q(y)^{1/2}d y \Big)
 \\
&\times \Big(1+O\big(\big|\int_{x_{\pm}}^x q(y)^{1/2}d y\big|^{-\rho_{1}}\big)\Big) 
\end{align*}
where $  \rho_{1}= \min\{\rho_{0}-1,1\}>0$ as $x\ri\pm\ii$. 
\end{proposition}
     
Uniform asymptotic formulas for $u_\pm (x )$ are given in the following assertion.
 
 \begin{theorem}\label{Airy}
 Let Assumptions~\ref{infi} and \ref{onewell}  hold. If $\pm x\geq \pm x_{\pm}$, then the   solutions $u_\pm  (x ) = u_\pm  (x; \lambda,\h )$ admit the representations 
  \begin{equation}
u_ \pm  (x )= |\xi_{\pm}'  (x)|^{-1/2} \A(\h^{- 2/3}\xi_{\pm} (x)) \big(1+\varepsilon_{\pm} (x;\h )\big) 
\label{eq:R8}\end{equation}
where the remainder satisfies  the estimate
 \begin{equation}
|\varepsilon_{\pm} (x; \lambda, \h )| \leq C\h (1+ |\xi_{\pm} (x)|)^{-\rho },  \q  
 \rho= 3 \min\{\rho_{0}  - 1, 1\} /2> 0.
\label{eq:R8x}\end{equation}
Let $x_{1}\in (x_{-}, x_{+})$. On the interval $  [ x_1, x_{+}  ]$ $($on the interval $  [ x_{-}   , x_1 ])$ the function $u_+$ $($the function $u_-)$ admits the representation 
\begin{equation}
u_\pm (x )= | \xi_{\pm}' (x)|^{-1/2} \A(\h^{-2/3}\xi_{\pm}(x)) + O(\h^{7/6} ( \h^{2/3} + |x-x_{\pm}|)^{-1/4}). 
\label{eq:R8-}\end{equation}
 \end{theorem}
 
 Away from the points $x_{\pm}$, we can replace the Airy function $\A (t)$ by its asymptotics  \e{eq:R6} or \e{eq:R7}. Indeed, in view of \e{eq:R4x}, we see that 
\begin{equation}
 |\xi_{\pm} (x) |\geq c |x-x_{\pm}|, \q c>0,
\label{eq:R4xx}\end{equation}
 and hence $\h^{-2/3} \xi_{+}(x)\to \pm \infty$ if $\h^{-2/3}   (x-x_{+})\to \pm \infty$ and  $\h^{-2/3} \xi_{-}(x)\to \pm \infty$ if $\h^{-2/3}   (x-x_{-})\to \mp \infty$.
This leads to the following result.
  
 \begin{corollary}\label{AiryX}
Suppose that $\delta_{\h}\h^{-2/3}\geq c>0$ $($in particular, $\delta_{\h} $ may be fixed$)$.
Then the  functions $u_\pm (x )$ have asymptotics
 \begin{equation}
u_\pm (x )= 2^{-1}  \pi^{1/2} \h^{1/6} q(x)^{-1/4} \exp\Big(\mp \h^{-1}\int_{x_{\pm}}^x q(y)^{1/2}dy \Big) \big(1+O(\h |\xi_{\pm} (x)|^{-3/2})\big)
\label{eq:R8simp}\end{equation}
as $\h\ri 0$ uniformly in $x\geq x_{+}+ \delta_{\h}$ for $u_+ (x)$ and in  $x\leq x_{-}-\delta_{\h}$ for $u_- (x) $. 
Let $x_{1} \in (x_{-},x_{+})$. Then
the functions $u_\pm (x )$    have asymptotics
  \begin{align}
u_\pm (x )  =    \pi^{1/2} \h^{1/6} |q(x)|^{-1/4} &  \sin \Big(\pm \h^{-1}\int^{x_{\pm}}_x | q(y)|^{1/2}d y  + \pi/4 \Big)
\nonumber\\
& + O(\h^{7/6} |x-x_{\pm}|^{-7/4})
\label{eq:R8-simp}\end{align}
as $\h\ri 0$ uniformly in   $x\in [  x_1, x_{+} -\delta_{\h}]$ for $u_+  (x) $ and uniformly in   $x\in [x_{-} + \delta_{\h} , x_1]$ for $u_-  (x)$.
 \end{corollary}

On the other hand, using estimates \e{eq:R4xx} and $|\A (t)| \leq C (1+ |t|)^{-1/4}$,  
we obtain   uniform in $\h$ estimates of the functions $u_\pm (x )$ in neighborhoods of the turning points.
 
 \begin{corollary}\label{AiryXX}
 For sufficiently small $|x-x_{\pm}|$, the estimate 
  \begin{equation}
|u_\pm (x ) |\leq C (1+ \h^{-2/3} | x-x_{\pm}|)^{-1/4}
\label{eq:turn}\end{equation}
holds with a constant $C$ which does not depend on $\h$.
 \end{corollary}

 We note that all asymptotic relations \e{eq:R8}, \e{eq:R8-}, \e{eq:R8simp} and \e{eq:R8-simp} can be differentiated with respect to $x$. In particular, we have   asymptotics
   \begin{align}
  u_\pm' (x )  =   \mp \pi^{1/2} \h^{-5/6} |q(x)|^{1/4} & \cos \Big(\pm \h^{-1}\int^{x_{\pm}}_x | q(y)|^{1/2}d y  + \pi/4 \Big)
\nonumber\\
& + O(\h^{1/6} |x-x_{\pm}|^{-7/4})
\label{eq:R8-sim}\end{align}
as $\h\ri 0$ uniformly in   $x\in [  x_1, x_{+} -\delta_{\h}]$ for $u_+  (x) $ and uniformly in   $x\in [x_{-} + \delta_{\h} , x_1]$ for $u_-  (x)$.
All these relations  can also be differentiated with respect to $\lambda$. For example, we have
   \begin{align*}
\partial  u_\pm (x ; \lambda,\h) / \partial&\lambda =   \pm   2^{-1} \pi^{1/2} \h^{-5/6} |q(x;\lambda)|^{-1/4}  
\int^{x_{\pm}}_x | q(y;\lambda)|^{-1/2}d y 
 \\
&\times\cos \Big(\pm \h^{-1}\int^{x_{\pm}}_x | q(y;\lambda)|^{1/2}d y  + \pi/4  \Big) + O(\h^{1/6} |x-x_{\pm}|^{-7/4})
 \end{align*}
as $\h\ri 0$ uniformly in   $x\in [  x_1, x_{+} -\delta_{\h}]$ for $u_+  (x) $ and uniformly in   $x\in [x_{-} + \delta_{\h} , x_1]$ for $u_-  (x)$.

 \medskip
 
 {\bf 2.5.}
 Let us now calculate the norm of the function $ u_\pm (x)$ in the space $L^{2}(x_{1}, \pm\infty)$ where $x_{1}\in (x_{-}, x_{+})$. Actually, we will obtain a more general result.
 
 \begin{proposition}\label{AiryN}
 Let a function $w(x)$ be differentiable on the interval $(x_{-}, x_{+})$ except a finite number of points. Suppose that $w(x)$and $w'(x)$ are locally bounded functions and that, for some $N$,
      \begin{equation}
|w(x)| \leq C q(x) \Bigl|  \int_{0}^{x } |q(y)| ^{1/2} dy  \Bigr|^N 
\label{eq:w}\end{equation}
if   $|x|$ is large.
 Then under Assumptions~\ref{infi} and \ref{onewell} we have the asymptotic relation
       \begin{equation}
\int_{x_1}^{\pm \infty} w(x) u_\pm^2   (x )  dx =   2^{-1 } \pi \h^{1/3}
\int_{x_1}^{x_{\pm}} w(x) (\lambda-v(x)) ^{-1/2} dx
 +O (\h^{2/3}  ). 
\label{eq:w1}\end{equation}
\end{proposition}
 
 \begin{pf}
 We will prove \e{eq:w1} for the sign $``+"$,  omit this index and add $\h$.
   Using asymptotics \e{eq:R8}, \e{eq:R8x} and the estimate
$\xi^\prime (x) \geq c> 0$,
 we see that
\[
\int_{x_{+}}^{x_{+}+1} w(x) u_{\h}^2  (x )  dx \leq C
\int_{x_{+}}^{\infty} \xi^\prime (x) \A^2 (\h^{-2/3}\xi  (x)) dx=C_{1}  \h^{2/3} .
\]
Similarly, using identity \e{eq:Req} and condition \e{eq:w}  we find that
\[
\int_{x_{+}+1}^\infty w(x) u_{\h}^2  (x )  dx \leq C
\int_{x_{+}+1}^{\infty} \xi^\prime (x) \xi^{3N/2 +1} (x) \A^2 (\h^{-2/3}\xi  (x))  dx= O (  \h^{\infty }) .
\]

Suppose that $ \d_{\h} \to 0$  as $ \h\to 0$ but $\d_{\h} \h^{-2/3}\geq c>0$.
The integral of $u_{\h}^{2} (x ) $ over $(x_{+}  -\delta_{\h},x_{+} )$   is   estimated by $C\d_{\h}$ because according to   \e{eq:turn} the functions $u_{\h}  (x ) $ are uniformly bounded in a neighborhood of the point $x_{+}$.  On the interval $(x_1 , x_{+}  -\delta_{\h})$, we have a relation
\begin{align}
\int_{x_1}^{x_{+} -\delta_{\h}} w(x)& u_{\h}^2  (x )  dx =
\pi \h^{1/3} \int_{x_1}^{x_{+} -\delta_{\h}}   w(x) | q(x)| ^{-1/2}  
\nonumber\\
\times
&\sin^2\Bigl(  \h^{-1} \int_{x}^{x_{+}} | q(y)| ^{1/2} dy +\pi/4 \Bigr)dx
  + O(\h^{4/3} \delta_{\h}^{-1}).
\label{eq:pm2b}\end{align}
Indeed, in view of asymptotics \e{eq:R8-simp} we have to show that the integrals
\[
\h^{7/3} \int_{x_1}^{x_{+} -\delta_{\h}}     | x-x_{+}| ^{-7/2} dx\q \mbox{and} \q
\h^{4/3} \int_{x_1}^{x_{+} -\delta_{\h}}     | q(x)| ^{-1/4}  | x-x_{+}| ^{-7/4} dx 
\]
are $O(\h^{4/3} \delta_{\h}^{-1})$. The first of them equals $C \h^{7/3} \delta_{\h}^{-5/2}$
which is $O(\h^{4/3} \delta_{\h}^{-1})$ because  $\h =O( \delta_{\h}^{ 3/2})$. To estimate the second integral,  we have to additionally take into account that
\begin{equation}
|q(x)| \geq c (x_{+} - x), \q c>0.
\label{eq:pmbq}\end{equation}

Next, we replace $sin^{2}(\cdot)$ in the right-hand side of \e{eq:pm2b} by $1/2$. Let us estimate the error.  Integrating by parts separately on every interval  where $w(x)$
is differentiable, we see that
\begin{align*}
\int_{x_1}^{x_{+} -\delta_{\h}} w(x) | q(x)|^{-1/2}& \exp\Bigl(2 i \h^{-1} \int_{x}^{x_{+}} | q(y)| ^{1/2} dy  \Bigr) dx
 \\
= - 2^{-1} i \h  w(x_{+}  -\delta_{\h})&  q(x_{+}  -\delta_{\h})^{-1}\exp\Bigl( 2 i \h^{-1} \int_{x_{+} -\delta_{\h}}^{x_{+}} |q(y)| ^{1/2} dy  \Bigr)
 \\
 + 2^{-1} i \h \int_{x_1}^{x_{+} -\delta_{\h}} &\big(  w'(x)   q(x) ^{-1}   - v'(x) q(x) ^{-2} w(x)  \big)\\\times
 & \exp\Bigl( 2 i \h^{-1} \int_{x}^{x_{+}} |q(y)| ^{1/2} dy  \Bigr) dx+ O(\h ).
\end{align*}
The right-hand side here is bounded by 
\[
C \h \Big(1+  | q(x_{+} -\delta_{\h})|^{-1}+\int_{x_1}^{x_{+} 
-\delta_{\h}} q(x) ^{-2}   dx\Big)  
\]
which in view of estimate \e{eq:pmbq}
   does not exceed $C \h \d_{\h}^{-1}$. Thus, it follows from \e{eq:pm2b} that
     \[
\int_{x_1}^{x_{+}-\d_{\h}} w(x) u_{\h}^{2}  (x )  dx = 2^{-1 } \pi \h^{1/3}
\int_{x_1}^{x_{+} -\d_{\h}}w(x)  (\lambda-v(x)) ^{-1/2} dx + O(\h^{4/3} \d_{\h}^{-1})   .
\]
Finally, making an error of order $O( \h^{1/3} \d_{\h}^{1/2})$, we can extend the integral in the right-hand side to the whole interval $(x_{1}, x_{+})$.
 Setting $\d_{\h}= \h^{2/3  }$   and putting   the results obtained together, we arrive at asymptotic relation \e{eq:w1}. 
    \end{pf}
    
    Of course \e{eq:w} is a very mild restriction.  It is satisfied  for $v(x)=w(x)$. It is also true for all    functions $v(x)$ if $w(x)$ if bounded by some power of $| x |$ at infinity and is even less restrictive if $v(x)\to \infty$ as $|x| \to \infty$. In particular, setting $w(x)=1$, we obtain
    
   \begin{corollary}\label{AiryNc}
The asymptotic relation holds:
       \begin{equation}
\int_{x_1}^{\pm\infty}  u_\pm^2  (x ) dx =   2^{-1 } \pi \h^{1/3}
\int_{x_1}^{x_{\pm}}(\lambda-v(x)) ^{-1/2} dx
 +O (\h^{2/3 }  ). 
\label{eq:pm3+}\end{equation}
\end{corollary}

\section{Proof of Theorem~\ref{Airy}}


 {\bf 3.1.}
 We will prove Theorem~\ref{Airy} for the sign $``+"$ and omit this index. On the contrary, we add the index $\h$ to emphasize the dependence on it of various objects.
 Let $x_{1}\in (x_{-},x_{+})$, $x\in (x_{1},\ii)$ and let the function $\xi (x)$ be defined by formulas \e{eq:R3}.
 According to Lemma~\ref{xi}   $\xi (x) \in (\xi_1, \infty)$   where $\xi_1= \xi(x_1)$ and    $x$ can be considered as a function of $\xi$ if $\xi  \in ( \xi_1, \infty)$.   
 
  Let us make   the change of variables $x\mapsto  \xi$ in  equation \e{eq:R1}  and set
\begin{equation}
u_{\h}  (x )= \xi^\prime(x)^{-1/2} f_{\h} (\h^{- 2/3}\xi(x))  .
\label{eq:R9}\end{equation}
Then using identity \e{eq:Req}, we obtain that
\begin{equation}
 - f_{\h}^{\prime \prime}(\h^{-2/3}\xi  ) +\h^{-2/3} \xi f_{\h} (\h^{-2/3}\xi  ) =\h^{4/3} r( \xi) f_{\h} (\h^{-2/3}\xi  ),
\label{eq:Ra1}\end{equation}
where
\begin{equation}
r( \xi)= p( x(\xi))    
\label{eq:rp}\end{equation} 
   and $  p ( x )$ is defined by formula \e{eq:Ra2}. 
  In view of   \e{eq:Ra2x}  we have the estimate 
\begin{equation}
  | r( \xi)|\leq C  (1+|\xi|)^{-1/2-\rho},\q    \rho= 3 \min\{\rho_{0}  - 1, 1\} /2> 0  .
\label{eq:Ra3}\end{equation}
 Setting in  \e{eq:Ra1} $t= \h^{-2/3}\xi$, we get the following intermediary result.

 \begin{lemma}\label{Airy1}
 Let $t= \h^{-2/3}\xi (x)$,  and let the functions $u_{\h }(x)$ and $f_{\h} ( t  )$ be related by formula \e{eq:R9}. Then
  equation \e{eq:R1} for $x\geq x_{1}$ is equivalent to the equation
 \begin{equation}
 - f_{\h}^{\prime \prime}( t  ) + t f_{\h} ( t  ) =
 R_{\h}   (t ) f_{\h} ( t ) \q \mbox{for} \q t\geq \xi_{1} \h^{-2/3}
\label{eq:Ra4}\end{equation}
where 
 \begin{equation} 
 R_{\h}   ( t ) = \h^{4/3} r(\h^{2/3}t ).
\label{eq:Ra5}\end{equation}
\end{lemma}

 \medskip
 
 {\bf 3.2.}
 Let us reduce differential equation \e{eq:Ra4} to a Volterra integral equation. Set 
\begin{equation}
K_{\h} (t,s )= - \pi 
\big(\A (t) \B (s) -\A (s) \B (t)\big)
 R_{\h}   (s ), \q s\geq t,
\label{eq:INT1}\end{equation}
and consider  the   equation
\begin{equation}
f_{\h} ( t  )= \A (t) +\int_{t}^\ii K_{\h} (t,s ) f_{\h} ( s  ) d s.
\label{eq:INT3}\end{equation}
Differentiating it twice, we see that its solution satisfies also  differential equation \e{eq:Ra4}. We will study equations \e{eq:Ra4} or \e{eq:INT3} separately for $t\geq 0$ and
$t\leq 0$.

 \begin{lemma}\label{Airy2}
 For $t\geq 0$, equation \e{eq:Ra4} has a solution $f_{\h} ( t  )$ such that
 \begin{equation}
f_{\h} ( t  )= \A (t) \big( 1+ \eta_{\h} (t )\big)
\label{eq:INT3f}\end{equation}
where
 \begin{equation}
|\eta_{\h} (t )|  \leq  C \h (1+ \h^{2/3} t)^{ -\rho },\q   \rho= 3 \min\{\rho_{0}  - 1, 1\} /2> 0.
\label{eq:INT3fX}\end{equation}
\end{lemma}

 \begin{pf} 
Making the multiplicative change of variables
\begin{equation}
f_{\h} (t )=   \A(t) g_{\h} (t ) 
\label{eq:Ra6}\end{equation}
and using \e{eq:R7cc}, we rewrite equation \e{eq:INT3} as
\begin{equation}
g_{\h} ( t  )= 1 -\int_{t}^\ii L_{\h} (t,s ) g_{\h} ( s  ) d s,
\label{eq:Int3}\end{equation}
where
\[
 L_{\h} (t,s ) = \A (t)^{-1}  K_{\h} (t,s ) \A (s)=   \int_{t}^s  \A^{-2}(\tau) d \tau
 \A^2(s)  R_{\h}   (s ), \q s\geq t. 
 \]
It follows from \e{eq:R6} that 
\[
\int_{t}^s  \A^{-2}(\tau) d \tau \leq C \exp (4s^{3/2}/3)  
\]
so that according to \e{eq:Ra3} and \e{eq:Ra5}  
\begin{equation}
| L_{\h} (t,s )| \leq C  \h^{4/3} s^{-1/2} (1+\h^{2/3}  s )^{-1/2-\rho},\q 0\leq t\leq s. 
\label{eq:Int2}\end{equation}
This estimate allows us to solve
 equation \e{eq:Int3} by iterations. In particular, the solution of \e{eq:Int3}   satisfies the estimate
\[
| g_{\h} ( t  ) -1 |Ê\leq C\int_{t}^\ii | L_{\h} (t,s ) | d s.
\]
Now   estimate \e{eq:INT3fX} on the remainder $\eta_{\h} (t )=g_{\h} (t )-1$ follows  again from \e{eq:Int2}.
 \end{pf}
 
 Putting together formulas \e{eq:R9} and  \e{eq:INT3f}, we obtain representation \e{eq:R8} with
$ \varepsilon_{\h}(x)= \eta_{\h} (\h^{-2/3} \xi(x) )$.
 Estimate \e{eq:INT3fX} implies estimate  \e{eq:R8x}. This leads to
   the assertion of Theorem~\ref{Airy} for $x\geq x_{+}$. In particular, for a fixed $\h$, we get Proposition~\ref{inf}.

  Next, we consider the case $t\leq 0$.
  
   \begin{lemma}\label{Airy3}
 For    $t \in  [   \h^{-2/3} \xi_{1}, 0]$,   the solution $f_{\h} ( t  )$ of equation \e{eq:Ra4}    satisfies the estimate
\begin{equation}
 | f_{\h} ( t  ) -  \A ( t  ) |Ê\leq  C  \h   (1+|t|)^{-1/4}.
\label{eq:Int4xx}\end{equation}
\end{lemma}

  \begin{pf}
   Let us rewrite equation \e{eq:INT3} as
\begin{equation}
f_{\h} ( t  )= f_{\h}^{(0)} ( t  ) +\int_{t}^0 K_{\h} (t,s ) f_{\h} ( s  ) d s,
\label{eq:INT3-}\end{equation}
where the new ``free" term
\begin{equation}
  f_{\h}^{(0)} ( t  ) =\A (t)+ f_{\h}^{(1)} ( t  ) ,\q
   f_{\h}^{(1)} ( t  ) =\int_{0}^\ii K_{\h} (t,s ) f_{\h} ( s  ) d s.
\label{eq:INT3fr}\end{equation}
It follows from \e{eq:R6}  and  \e{eq:R6B}  that
\[
 \A^2 (t)+  \A (t)  \B(t) \leq C (1+t)^{-1/2}, \q t\geq 0,
\]
and from \e{eq:R7}  and  \e{eq:R7B}  that
\begin{equation}
 |\A (t)| +  | \B(t)|\leq C (1+|t|)^{-1/4}, \q t\leq 0.
\label{eq:INT4-}\end{equation}
Therefore using \e{eq:Ra3}, \e{eq:Ra5} and \e{eq:INT1}, we find that 
\begin{align}
|f_{\h}^{(1)} ( t  )| &\leq  C \Big( |\A (t)| \int_{0}^\ii  \A (s) \B (s) | R_{\h} (s )| d s
+ |\B (t)| \int_{0}^\ii    \A^2 (s) | R_{\h} (s )| d s \Big)
\nonumber\\
&\leq 
 C_{1} (1+|t|)^{-1/4} \h^{4/3}\int_{0}^\ii   s^{-1/2} (1+\h^{2/3}s )^{-1/2-\rho}   ds
\nonumber\\
&\leq 
 C_{2} (1+|t|)^{-1/4} \h.
\label{eq:INT3fr1}\end{align}

Let us now consider  equation  \e{eq:INT3-}.
By virtue of estimates \e{eq:Ra3} and \e{eq:INT4-} its kernel satisfies the bound
\begin{equation}
|K_{\h}  (t,s )|\leq C \h^{4/3}
 (1+|t|)^{-1/4} (1+|s|)^{-1/4} r(\h^{2/3} s), \q t\leq s\leq 0 , 
\label{eq:Int2+}\end{equation}
where the function $r(\h^{2/3} s)$ can be estimated by a constant.
Thus, solving   \e{eq:INT3-} again by iterations, we obtain the estimate
\begin{align}
| f_{\h} ( t  ) -  f_{\h}^{(0)} ( t  ) |Ê&\leq  C_{1}\int_{t}^0 | K_{\h} (t,s ) | (1+|s|)^{-1/4} d s
\leq C_{2} \h^{4/3}  \int_{t}^0 (1+|s|)^{-3/4} d s 
\nonumber\\
&\leq   C_{3} \h^{4/3}  (1+|t|)^{1/4}   .
\label{eq:Int4+}\end{align}
If $t \in  [   \h^{-2/3} \xi_{1}, 0]$, then combining definition   \e{eq:INT3fr} with estimates  \e{eq:INT3fr1} and \e{eq:Int4+}, we get estimate \e{eq:Int4xx}.
   \end{pf} 
  
 In view of formula \e{eq:R9}, this lemma yields the result of Theorem~\ref{Airy} for 
 $x\in [x_{1} , x_{+}]$.
 
 Differentiating integral equation  \e{eq:INT3} with respect to $t$, we obtain
 asymptotic relations for        $f_{\h}' (t)$ and then for        $u_{\h}' (x)$. 
This concludes the proof of Theorem~\ref{Airy}.

\section{Semiclassical asymptotics of eigenfunctions}


 {\bf 4.1.} 
  Let $ \lambda =\lambda(\h)$ be an eigenvalue of the   Schr\"odinger operator $H_{\h} = -\h^2 d^2/ dx^2 + v(x)$ from a neighborhood of a non-critical point $\lambda_{0}$.  Then the solutions $u_\pm (x) $ are proportional:
        \begin{equation}
u_- (x ; \lambda,\h) = a (\lambda ,\h) u_+ (x ; \lambda,\h) . 
\label{eq:E1}\end{equation}
Choose     an arbitrary interior point $x$ of the interval $(x_{-} (\lambda ),x_{+}(\lambda ))$. To  calculate   the Wronskian of $u_+  (x)$ and $u_-  (x)$,  we  use asymptotic relations \e{eq:R8-simp} and \e{eq:R8-sim}. 
 Setting
     \begin{equation}
\varphi_{\pm}  (x; \lambda )=\pm \int^{x_{\pm}(\lambda )}_x (\lambda  - v(y))^{1/2}dy   , 
\q x\in (x_{-} (\lambda ),x_{+}(\lambda )), 
\label{eq:E1ph}\end{equation}
we find that
  \begin{align}
w(\lambda,\h) = u_ +  (x;\lambda,\h)   u&'_ -  (x;\lambda,\h)  -  u_ - (x;\lambda,\h)   u'_+  (x;\lambda,\h) 
\nonumber \\
=
\pi \h^{-2/3}  \Big( \sin &(\h^{-1}\varphi_{+} (x; \lambda )+\pi/4 )
  \cos (\h^{-1}\varphi_{-}  (x; \lambda )+\pi/4 )
 \nonumber \\
   + \cos (\h^{-1}\varphi&_{+}  (x; \lambda  )+\pi/4 )
  \sin (\h^{-1}\varphi_{-}  (x; \lambda  +\pi/4)\Big) + O(\h^{1/3})
\nonumber   \\
=
\pi  & \h^{-2/3}   \sin (\h^{-1}\Phi (\lambda ) +\pi/2)+ O(\h^{1/3})
\label{eq:Ewr}  \end{align}
where $\Phi (\lambda )=\varphi_{+} (x; \lambda )+\varphi_{-}  (x; \lambda )$ so that
   \begin{equation}
\Phi (\lambda )  =\int^{x_{+}(\lambda )}_{x_{-}(\lambda )}  (\lambda - v(y))^{1/2}d y.
\label{eq:Phi}\end{equation}
Since $w(\lambda,\h) =0$, we see that 
\[
\sin (\h^{-1}\Phi (\lambda ) +\pi/2)= O(\h )
\]
 and hence
 \begin{equation}
\int^{x_{+}(\lambda )}_{x_{-}(\lambda )}  (\lambda  - v(x))^{1/2}d x= \pi ( n  +1/2) \h  + O(\h^2)
\label{eq:B-S}\end{equation}
for some integer number $n=n (\lambda, \h ) $. This gives us the famous Bohr-Sommerfeld quantization condition.

Suppose now that a number $\pi( n+1/2)  \h$ belongs to a   neighborhood  of $\lambda_{0}$. Let us check that there exists an eigenvalue $\lambda_{n}(\h)$ of the operator $H_{\h}$ satisfying the estimate 
  \begin{equation}
| \Phi (\lambda_{n}(\h)) -\pi( n+1/2)  \h|\leq C \h^2.
\label{eq:con}\end{equation}
Since $u_{\pm}\in L^2 ({\Bbb R}_{\pm})$, it suffices to show that  $w(\lambda,\h) =0$ for some $\lambda=\lambda_{n}(\h)$  satisfying   estimate \e{eq:con}.
Using the equality $\lambda - v(x_\pm(\lambda ))=0$, we find that 
  \begin{equation}
\Phi'(\lambda ) = 2^{-1}\int^{x_{+}(\lambda )}_{x_{-}(\lambda )}  (\lambda - v(y))^{-1/2}
d y>  0.
\label{eq:co}\end{equation}
 Hence $\Phi$ is a one-to-one  mapping of a   neighborhood  of $\lambda_{0}$ on a   neighborhood  of $\mu_{0}=\Phi (\lambda_{0})$. Set $\mu=\Phi (\lambda)$ and
\begin{equation}
\epsilon(\mu,\h)=  \pi^{-1} \h^{2/3} w(\Phi^{-1}(\mu),\h) -\sin (\h^{-1} \mu + \pi/2 ) .
\label{eq:con1}\end{equation}
In view of  \e{eq:Ewr} this function satisfies the estimate
$|\epsilon(\mu,\h)| \leq C\h$
with a constant $C$ which does not depend on $\h$ and
  $\mu$ from a   neighborhood  of $\mu_{0}$. We have to show that the equation
\[
\sin (\h^{-1} \mu+\pi/2) + \epsilon(\mu,\h)=0
\]
has a solution $\mu_{n}(\h)$ obeying the estimate
\[
|\mu_{n}(\h)-\pi( n+1/2)  \h|\leq C \h^2.
\]
Setting $s =\h^{-1} \mu+\pi/2$, we see that this assertion is equivalent to the existence of a solution $s=s_{n}(\h)$ of the equation
\begin{equation}
\sin s+ \epsilon(\h (s-\pi/2),\h)=0
\label{eq:con2}\end{equation}
obeying the estimate 
\begin{equation}
|s_{n}(\h)-\pi ( n +1) |\leq C \h .
\label{eq:con3}\end{equation}
The last fact is obvious because $\epsilon(\h (s-\pi/2),\h)=O(\h)$.  

Next, we will show that for every $n$ there is only one eigenvalue of the operator $H_{\h}$ satisfying \e{eq:con}. To that end, we have to check that equation \e{eq:con2} cannot have two solutions satisfying \e{eq:con3}. Supposing the contrary, we find  a point $\tilde{s}= \tilde{s}_{n}(\h)$ such that 
\begin{equation}
\cos \tilde{s}= -\h  \frac{\partial\epsilon}{\partial\mu}(\h (\tilde{s}-\pi/2),\h) 
\label{eq:con4}\end{equation}
and $\tilde{s}_{n}(\h) = \pi ( n +1)  + O( \h)$. Observe that relation \e{eq:Ewr}  can be differentiated in $\lambda$ which yields
\[
\partial w(\lambda,\h) / d\lambda=
\pi   \Phi'(\lambda)\h^{-5/3}   \cos (\h^{-1}\Phi (\lambda ) +\pi/2)+ O(\h^{-2/3}).
\]
It follows that function \e{eq:con1} obeys the estimate $\partial\varepsilon(\mu,\h)/d\mu= O(1)$. Thus, the right-hand side of equation \e{eq:con4} is $O(\h)$ while its left-hand side tends to $(-1)^{n+1}$ as $\h\to 0$. 

Finally, plugging asymptotics \e{eq:R8-simp} and \e{eq:R8-sim} into \e{eq:E1}, we find that
  \begin{align}
   \sin (\h^{-1} \varphi_{-}  (x; \lambda )+ \pi/4) &+O(\h)  
 \nonumber\\
 =
 & a(\lambda, \h )   \big(\sin (\h^{-1} \varphi_{+}  (x; \lambda )+ \pi/4) +O(\h)\big)
\label{eq:pm+}\end{align}
and
  \begin{align}
  \cos(\h^{-1} \varphi_{-} (x; \lambda )+ \pi/4) &+O(\h) 
  \nonumber\\
  = 
  & - a(\lambda, \h )  \big(\cos (\h^{-1} \varphi_{+}  (x; \lambda )+ \pi/4) +O(\h)\big).
\label{eq:pm-}\end{align}
Together,  these two relations imply that $|a (\lambda, \h)|=1 + O(\h)$.  Moreover, since
\[
\varphi_{+} (x; \lambda )+\varphi_{-}  (x; \lambda ) = \pi ( n  +1/2) \h  + O(\h^2),
\]
it follows from \e{eq:pm+} and \e{eq:pm-}   that
   \begin{equation}
a (\lambda, \h) = (- 1)^{n  }+ O(\h).
\label{eq:pm}\end{equation}

 Thus, we have obtained the following result.
 
  \begin{theorem}\label{pm}
  Let Assumptions~\ref{infi} and \ref{onewell} hold for  a point $\lambda_{0}$.
Suppose that an eigenvalue $ \lambda=\lambda(\h) $ of the operator $H_{\h}$ belongs to a neighborhood of  $\lambda_{0}$.
 Then necessarily condition  \e{eq:B-S}   is satisfied with some integer  number $n= n (\lambda ,\h) $. Conversely, for every $n$ such that $\pi( n+1/2)  \h$ belongs to a   neighborhood  of $\Phi (\lambda_{0})$, there exists an eigenvalue $\lambda_{n}(\h)$ of the operator $H_{\h}$ satisfying   estimate 
\e{eq:con}
with a constant $C$ not depending on $n$ and $\h$. Such an eigenvalue $\lambda_{n}(\h)$ is unique. Moreover, the coefficient $ a (\lambda, \h)  $ in \e{eq:E1} has asymptotics  \e{eq:pm} where $n $ is the same number as in \e{eq:B-S}.
   \end{theorem}

  \begin{corollary}\label{pmc}
 Let an interval $(a_{1},a_{2})$  belong to a neighborhood of a point $\lambda_{0}$ satisfying Assumptions~\ref{infi} and \ref{onewell}.  Then the total number $N_{\h}$ of eigenvalues  of the operator $H_{\h}$ in this interval equals
     \begin{equation}
  N_{\h}=\pi^{-1}( \Phi (a_{2}) -\Phi (a_1) ) \h^{-1}+ \epsilon (\h)
\label{eq:Weyl}\end{equation}
where $|\epsilon (\h)|\leq 1$   for sufficiently small $\h$.
   \end{corollary}

  \begin{pf}
  According to Theorem~\ref{pm} there is exactly one eigenvalue of the operator $H_{\h}$ in a neighborhood of size $C\h^2$ of every point $\Phi^{-1} (\pi (n+1/2)\h)$. These neighborhoods have empty intersections for sufficiently small $ \h$.
  Thus, $N_{\h}$ equals the number of points $ \pi (n+1/2)\h$ lying in the interval $(\Phi (a_1), \Phi (a_2))$. Clearly, this number equals the right-hand side of \e{eq:Weyl}.
      \end{pf}
      
       \begin{remark}\label{pmd}
       Suppose that Assumptions~\ref{infi} and \ref{onewell} hold true for all $\lambda\in [a_{1},a_{2}]$. Then remainders in different asymptotic formulas of this paper can be estimated uniformly  in  $\lambda\in [a_{1},a_{2}]$.  Formula \e{eq:Weyl} also remains true for such $(a_{1},a_{2})$. 
   \end{remark}

      Note that definition \e{eq:Phi} can be rewritten as
       \begin{equation}
      \Phi(\lambda)= 2^{-1} \int\int_{  p^2+ v(x) \leq \lambda}dpdx.
   \label{eq:Weyl1}\end{equation}
      Indeed, integrating in the right-hand side first over $p$ we obtain the right-hand side
      of \e{eq:Phi}. It follows that the asymptotic coefficient in \e{eq:Weyl} is the volume of a part of the phase space:
            \[
      \Phi (a_{2}) -\Phi (a_1) = 2^{-1} \mathrm{mes} \{(x,p)\in{\Bbb R}^2:   a_{1}\leq p^2+ v(x) \leq a_{2}\} .
      \]
        Thus, relation  \e{eq:Weyl} is the semiclassical Weyl formula with a strong estimate of the remainder.

 \medskip
 
 {\bf 4.2.}
Let us denote by $\psi (x)=\psi(x;\lambda,\h)$  the   eigenfunction of the operator $H_{\h}$ corresponding to its eigenvalue  $ \lambda $.    We suppose that $\psi =\overline{\psi } \in L^2 ({\Bbb R}_{+}) $ and $\|\psi \|=1$ which fixes $\psi $ up to a sign. Clearly,
   \begin{equation}
\psi  (x) =c_{\pm}   u_\pm (x ), \q c_{\pm}  =c_{\pm} (\lambda, \h ), 
\label{eq:E3}\end{equation}
where according to \e{eq:pm}
 \[
|c_{+} (\lambda,\h ) | = |c_{-}(\lambda, \h) | (1+ O(\h )).
\]
Therefore it follows from Corollary~\ref{AiryNc} that
   \begin{equation}
 |c_{\pm} (\lambda, \h ) | = 2^{1/2} \pi^{-1/2} \h^{-1/6}
\Big(\int_{x_{-}(\lambda )}^{x_{+}(\lambda )} (\lambda -v(x))^{-1/2} dx \Big)^{-1/2}+O(\h^{1/6   }),  
\label{eq:E4}\end{equation}
which in view of Theorem~\ref{Airy} yields the following result.
 
  \begin{theorem}\label{norm}
 Under the assumptions of Theorem~\ref{pm},
let us denote by $\psi (\lambda, \h)$ the real normalized eigenfunction $($defined up to a sign$)$ of the operator $H_{\h}$ corresponding to its eigenvalue  $ \lambda =\lambda(\h)$.  Let $x_{1}$ be an arbitrary point from the interval $(x_{- } (\lambda ), x_{+ } (\lambda ))$. Then,  for $x\in (x_{1 } ,\ii)$,  asymptotics of $\psi (x; \lambda, \h)$ as $\h\ri 0$  is given by formulas \e{eq:E3}, \e{eq:E4} for the sign $``+"$ and asymptotic relations of Theorem~\ref{Airy} for the function $u_+  (x ;\lambda, \h ) $. Similarly,   for $x\in (-\ii, x_{1 }  )$,  asymptotics of $\psi (x;\lambda, \h)$ as $\h\ri 0$  is given by formulas \e{eq:E3}, \e{eq:E4}  for the sign $``-"$ and asymptotic relations of Theorem~\ref{Airy}  for the function $u_-  (x ;\lambda, \h ) $.  In neighborhoods of the turning points,     the estimate 
  \begin{equation}
|\psi (x; \lambda, \h) |\leq C (  \h^{2/3} + | x-x_{\pm}|)^{-1/4}
\label{eq:turn1}\end{equation}
holds with a constant $C$ which does not depend on $\h$.
\end{theorem}

   In view of formula \eqref{eq:w1}, this result can be supplemented by the following
   
    \begin{proposition}\label{normw}
    Let a function $w$ satisfy the assumptions of Proposition~\ref{AiryN}. Then
  under the assumptions of Theorem~\ref{norm}  we have  
\begin{align}
  \int_{-\ii}^\ii  w(x) \psi ^2(x; \lambda,\h) dx  
 =&
  \int _{x_{-}(\lambda ) }^{x_{+}(\lambda ) } w(x) ( \lambda  -v(x))^{-1/2} dx 
\nonumber\\
&\times \Bigl( \int_{x_{-}(\lambda ) }^{x_{+}(\lambda ) }   (\lambda  -v(x))^{-1/2} dx \Bigr)^{-1}
+ O(\h^{1/3 }).
\label{eq:mainOld1}\end{align}
   \end{proposition}
   
     In particular, this relation applies to the potential energy 
     \[
 V(\lambda,\h)=     \int_{-\ii}^\ii  v(x) \psi ^2(x; \lambda,\h) dx
      \]     
     and by virtue of the energy conservation $K(\lambda,\h)+V(\lambda,\h)=\lambda$, we also obtain     the asymptotics   of the kinetic energy.
   
   \begin{corollary}\label{KIN}
  Under the assumptions of Theorem~\ref{norm}   asymptotic relation   \e{eq:mainkb} holds with the leading term $K_{cl} (\lambda)$ given by \e{eq:mainOld}.
  \end{corollary}
  
  Since $K_{cl} (\lambda)>0$, for small $\h$  the kinetic energy $K(\lambda,\h)\geq c>0$ or, equivalently, 
  the potential energy $V(\lambda,  \h)\leq \lambda-c$. This implies that the eigenfunctions $\psi (x; \lambda,\h)$ are not too strongly localized in neighborhoods of the turning points 
  $x_{\pm } (\lambda )$. In view of estimate \e{eq:turn1}, this statement can be reinforced.
  
    \begin{proposition}\label{loc}
  Let the assumptions of Theorem~\ref{pm} hold, and let $\| \psi (  \lambda, \h)\|=1$. Then
\[
  \int_{x_{\pm} (\lambda ) -\d}^{x_{\pm} (\lambda )+\d}  \psi ^2(x; \lambda, \h) dx  
\leq C \d^{1/2}
\]
where the constant $C$   does not depend on $\h$.
   \end{proposition}
    
  \begin{remark}\label{unbound}
  It follows from asymptotics   \e{eq:R8} and   \e{eq:E4} that
   $$
  | \psi (x_\pm( \lambda); \lambda , \h) | = \alpha_{\pm}( \lambda )\h^{-1/6}
   (1+O (\h^{1/3 })),  
   $$
   where the coefficient
   \[
   \alpha_{\pm}( \lambda )= 2^{1/2} \pi^{-1/2}  
\Big(\int_{x_{-}(\lambda )}^{x_{+}(\lambda )} (\lambda -v(x))^{-1/2} dx \Big)^{-1/2} |v'(x_{\pm} (\lambda))|^{-1/6} \A (0) \neq 0.
    \] 
  This contradicts the assertion of Theorem~7.1 of \cite{Simon}  that normalized eigenfunctions are uniformly bounded in neighborhoods of   turning points.
   \end{remark}

\medskip
 
 {\bf 4.3.}
    Recall that  a classical particle (of mass $m$ and energy $\lambda$) moves periodically (see, e.g.,  \cite{LLcl}) in a potential well bounded by the points $x_{-}=x_{-}(\lambda)$ and $x_{+} = x_{+}(\lambda)$ such that $v(x_{\pm}) =\lambda$. Let us check that the asymptotic coefficient 
 $K_{cl} $ in   \e{eq:mainkb}  
  coincides with the averaged over the period $T=T(\lambda)$ value 
  \[
  K_{av} =T ^{-1}\int_{0}^{T } K(t)dt
  \]
    of the classical  kinetic energy 
 \[
K(t)=mx^\prime(t)^2/2= \lambda-v(x(t)).
\]
Since
\[
dt=x'(t)^{-1} dx= (m/2)^{1/2}(\lambda-v(x)) ^{-1/2} dx,
\]
  the period  is given by the formula 
\[
T = 2 \int_{x_{-}}^{x_{+}} \frac{dt}{dx} dx =(2m)^{1/2}  \int_{x_{-}}^{x_{+}}  (\lambda-v(x))^{-1/2} dx  
\]
and
\[
K_{av} =T ^{-1} \int_0^{T}  (\lambda-v(x(t)))  dt
= 2  (m/2)^{1/2}T ^{-1} \int_{x_{-} }^{x_{+}}  (\lambda-v(x))^{1/2}  dx. 
\]
Putting these two relations  together, we obtain for $K_{av} $ the same expression   \e{eq:mainOld} as for $K_{cl} $. This proves the equality
  \[
K_{cl}(\lambda)=K_{av}(\lambda ).
\]

We also note that    
\[
K_{cl}(\lambda)=   \big( 2 d \ln\Phi(\lambda)/ d\lambda\big)^{-1}
\]
where the function $\Phi(\lambda)$ is defined by formulas  \e{eq:Phi} or, equivalently,
 \e{eq:Weyl1}. For the proof of this equality, it suffices to plug representation \e{eq:co} for  the function $\Phi'(\lambda)$   into formula \e{eq:mainOld}.
  
\section{Discontinuous potentials}

 {\bf 5.1.}
 Away from the turning points,
 assumptions on $v(x)$ can be somewhat relaxed.  Consider,  for example,  an interval $(x_{-}+ \d, x_{+} - \d )$ where $\d>0$.  There, it suffices to require    that $v\in C^1$ and that $v'$ be absolutely continuous so that $v'' \in L^1$ (instead of 
 $v\in C^2$). In this case the function $r(\xi)$ defined by formulas \e{eq:Ra2} and \e{eq:rp} belongs to $L^1$ only so that the factor $r(\h^{2/3}s)$ in the right-hand side of estimate \e{eq:Int2+} cannot be neglected. Therefore (cf.  \e{eq:Int4+}) we have the estimate
 \begin{align*}
 \int_t^0 |K_{\h} (t,s)| (1+|s|)^{-1/4}ds
& \leq C \h^{4/3} (1+|t|)^{-1/4} \int_t^0 |r(\h^{2/3}s)|  ds
\\
&  \leq C_{1} \h^{2/3} (1+|t|)^{-1/4} \int_{\xi (x_{-}+ \d)}^0 | r( s) |  ds. 
 \end{align*}
 It follows that instead of  \e{eq:Int4xx} we have a slightly weaker estimate with $\h^{2/3}$ in place of $\h $ in the right-hand side. All other estimates remain unchanged.  Thus, Theorem~\ref{Airy}  is true   with a little bit weaker estimates of remainders in asymptotic formulas  for $u_{\pm}(x;\lambda,\h)$ inside the interval $(x_{-}+ \d, x_{+} - \d )$.   Repeating the arguments of Section~4, we get the following result.
 
  \begin{proposition}\label{Ddisc}
Under the assumptions above, all results of Theorem~\ref{pm} $($and of Corollary~\ref{pmc}$)$ about eigenvalues   of the operators $H_{\h}$ remain true with the remainders $O (\h^{5/3}) $  in \e{eq:B-S},  $O (\h^{2/3}) $   in  \e{eq:pm}  and $C\h^{5/3}$ in the right-hand side of  \e{eq:con}.   Theorem~\ref{norm} about corresponding eigenfunctions      remains also true with the remainders $O (\h^{5/6} |x-x_{\pm}|^{-7/4})$ in \e{eq:R8-simp},  $O (\h^{-1/6} |x-x_{\pm}|^{-7/4})$ in  \e{eq:R8-sim} and $O (\h^{1/18})  $ in  \e{eq:E4}.
\end{proposition}

\medskip
 
 {\bf  5.2.}
 Our goal in this subsection is to extend the results of Section~4 to functions $v(x)$ with a singular point $x_{0}$ inside a potential well.
 
We suppose that  Assumption~\ref{infi} holds everywhere except a point $x_{0}$  and that  Assumption~\ref{onewell} holds   for some $\lambda_{0}$ such that $x_{0}$ is an interior point of the interval $(x_{-}(\lambda_{0}), x_{+}(\lambda_{0}))$. We   assume that $v(x)$   has finite limits at   $x_{0}$  but the left and right limits might be different. Finally, we require that $v' \in L^2 (x_{0} , x_{0} \pm \d)$ and $v'' \in L^1 (x_{0} , x_{0} \pm\d)$  for some $\d>0$. 
 
 Now    we can construct solutions $u_+(x)$ and $u_-(x)$ of equation \e{eq:I2} on the intervals $(x_{0}, \infty)$ and $(-\infty, x_{0})$, respectively.  Define, as usual,  the function $r(\xi)$   by formulas \e{eq:Ra2} and \e{eq:rp}. Since $r\in L^1(x_{0} , x_{0} \pm \d)$, 
   the limits 
 $ u_\pm (x_{0} \pm 0)$ and $ u_\pm' (x_{0} \pm 0)$ 
 exist, and we can use formulas  \e{eq:R8-simp} and  \e{eq:R8-sim} for these limits
 (with slightly weaker estimates of the remainders -- see subs.~5.1).
 It follows that  the Wronskian $w(\lambda,\h)$  of $u_+$ and $u_-$ calculated at the point $x_{0}$ is   given by the expression  (cf. \e{eq:Ewr})
   \begin{align*}
   \pi \h^{-2/3}  \Big(   
     p(x_{0}, \lambda ) \sin (\h^{-1}\varphi&_{+}  (x_{0}; \lambda )+\pi/4 )
  \cos (\h^{-1}\varphi_{-}  (x_{0}; \lambda )+\pi/4 )
  \\
  +
  p(x_{0}, \lambda )^{-1} 
   \cos & (\h^{-1}\varphi_{+}  (x_{0}; \lambda )+\pi/4 )
  \sin (\h^{-1}\varphi_{-}  (x_{0}; \lambda ) +\pi/4)
   \Big)
   + O(1) 
  \end{align*}
where 
  \begin{equation}
 p(x_{0}, \lambda )=   
  \big(\lambda -v(x_{0}-0) \big)^{1/4}    \big(\lambda -v(x_{0}+ 0) \big)^{-1/4} .
 \label{eq:PP} \end{equation} 
 
  Let an eigenvalue $\lambda$ of the operator $H_{\h}$ be close to $\lambda_{0}$. Since
  $w(\lambda,\h)=0$, we see that  
  \begin{align}
p(x_{0}, \lambda ) \sin (\h^{-1}\varphi_{+}  (x_{0}; \lambda )+\pi/4 )
  \cos (\h^{-1}\varphi_{-}  (x_{0}; \lambda )+\pi/4 )&
  \nonumber\\
   +   p(x_{0}, \lambda )^{-1} \sin (\h^{-1}\varphi_{+}  (x_{0}; \lambda )+\pi/4 )
  \cos (\h^{-1}\varphi_{-}  (x_{0}; \lambda )+\pi/4 )&
   = O (\h^{2/3}).
\label{eq:EwrDC}  \end{align}
 Formula \e{eq:EwrDC}  yields a generalization of the Bohr-Sommerfeld quantization condition \e{eq:B-S} and reduces to it if  $v(x_{0}+0)= v(x_{0}-0)$.

 Let the coefficient $a(\lambda, \h)$ be defined by equality \e{eq:E1}.  To calculate 
 $|a(\lambda, \h)|  $, we use again  relations \e{eq:pm+} and \e{eq:pm-}. However, additional factors $| q(x_{0}-0)|^{-1/4}$ and $| q(x_{0}-0)|^{1/4}$ appear now in their left-hand sides. Similarly, additional factors $| q(x_{0}+0)|^{-1/4}$ and $| q(x_{0}+0)|^{1/4}$ appear  in their right-hand sides. This implies  that
\begin{align}
a^2(\lambda, \h)= & p^2(x_{0}, \lambda )  \cos^2 (\h^{-1}\varphi_-  (x_{0}; \lambda )+\pi/4 ) \nonumber\\
&+ p^{-2}(x_{0}, \lambda )  \sin^2 (\h^{-1}\varphi_-  (x_{0}; \lambda )+\pi/4 ) + O (\h^{2/3})
\nonumber\\
=&\big( p^2(x_{0}, \lambda )  \sin^2 (\h^{-1}\varphi_+  (x_{0}; \lambda )+\pi/4 )
\nonumber\\ & + p^{-2}(x_{0}, \lambda )  \cos^2 (\h^{-1}\varphi_+  (x_{0}; \lambda )+\pi/4 )\big)^{-1}  + O (\h^{2/3}). 
 \label{eq:aa}  \end{align}
 As before, using formula \e{eq:pm3+} (where $x_{1}= x_{0}$) and the normalization condition $\| \psi \|=1$, we obtain explicit expressions for the absolute values of constants $c_{\pm} (\lambda, \h)$ in \e{eq:E3}: 
  \begin{align}
 |c_{+}(\lambda, \h) | =  2^{1/2}& \pi^{-1/2} \h^{-1/6}
\Big(    \int_{x_0}^{x_{+}(\lambda )} (\lambda -v(x))^{-1/2} dx 
\nonumber\\
+  &a^{-2}(\lambda, \h)\int_{x_{-}(\lambda )}^{x_{0} } (\lambda -v(x))^{-1/2} dx \Big)^{-1/2}+O(\h^{1/18   }) \label{eq:E4di}\end{align}
and 
 \begin{align}
 |c_{-}(\lambda, \h) | =  2^{1/2}& \pi^{-1/2} \h^{-1/6}
\Big( a^2(\lambda, \h)  \int_{x_0}^{x_{+}(\lambda )} (\lambda -v(x))^{-1/2} dx 
\nonumber\\
+  &\int_{x_{-}(\lambda )}^{x_{0} } (\lambda -v(x))^{-1/2} dx \Big)^{-1/2}+O(\h^{1/18   }). \label{eq:disc1}\end{align}

 Thus, Theorems~\ref{pm} and \ref{norm} can be supplemented by the following result.

 \begin{theorem}\label{disc}
Under the assumptions above, let an eigenvalue $\lambda=\lambda(\h)$ of the operator $H_{\h}$ belong to a neighborhood of $\lambda_{0}$.  Then necessarily   condition \e{eq:EwrDC}
  is satisfied with    the numbers $\varphi_{\pm}(x_{0};\lambda)$ and $p(x_{0},\lambda)$ defined by  \e{eq:E1ph}  and \e{eq:PP}, respectively. All assertions  $($for $x_{1}=x_{0})$ of Theorem~\ref{norm}  about the corresponding normalized eigenfunction $\psi(x;\lambda,\h)$ are true with the constants $c_{\pm} (\lambda, \h)$ whose absolute values are determined by   formulas
   \e{eq:aa},   \e{eq:E4di} and \e{eq:disc1}.
   \end{theorem}
   
     \begin{remark}\label{discR}
  If the functions $v'(x)$ and $v''(x)$  are bounded in a neighborhood of the point $x_{0}$, then even in the case   $v(x_{0}+0)\neq v(x_{0}-0)$ estimates of all remainders are the same as in Section~4. Thus, we have $O(\h)$ in \e{eq:EwrDC} and $O(\h^{1/6})$ in   \e{eq:aa} -- \e{eq:disc1}.
  \end{remark}

 \begin{remark}\label{discC}
 Let under the assumptions above $v(x_{0}+0)= v(x_{0}-0)$. Then all conclusions of Proposition~\ref{Ddisc}   remain true although the function $v'(x)$ is not required to be continuous at the point $x_{0}$.   In particular, we see that   jumps of derivatives of the function $v(x)$   at the point $x_{0}$ are inessential.
   \end{remark}

\medskip
 
 {\bf  5.3.}
 Let us consider an explicit example: 
 \begin{equation}
v(x)=a_{+} +  v_+ x^{\alpha_+}\;\; \mathrm{for}\;\; x>0 \q 
\mathrm{and} \;\; v(x)= a_{-} +v_-|x|^{\alpha_-} \;\; \mathrm{for}\q x<0, 
\label{eq:pm1}\end{equation}
where $v_{\pm}>0$ and $\alpha_{\pm}>0$.  Then all $\lambda>\max\{a_{+}, a_{-}\}$ are non-critical, the equation $v(x)=\lambda$ has two solutions $x_{+}> 0$, $x_{-}<0$ and $(x_{-}, x_{+})$ is a potential well. The point $x_{0}=0$ might be singular and 
$ p(0, \lambda)= (\lambda-a_{-})^{1/4} (\lambda-a_+)^{-1/4}$.

For potentials \e{eq:pm1}, the integrals in formulas \e{eq:Phi} and \e{eq:E4} can be calculated in terms of the beta function $\mathrm{B}$. 
Observe that $x_{+}= (\lambda v^{-1})^{1/ \alpha}$ if $v(x)= v  x^{\alpha }$ for $x>0$.  For the integrals over $(0,x_{+})$, we have  
 \begin{equation}
 \int _{0}^{x_{+}}  ( \lambda-vx^\alpha)^{1/2} dx=\lambda^{1/2} (\lambda/v)^{1/ \alpha} 
 \alpha^{-1}\mathrm{B}(3/2,1/\alpha)
 \label{eq:IN1}\end{equation}
 and 
 \begin{equation}
 \int _{0}^{x_{+}}  ( \lambda-vx^\alpha)^{-1/2} dx=\lambda^{-1/2} (\lambda/v)^{1/ \alpha} 
 \alpha^{-1}\mathrm{B}(1/2,1/\alpha).
\label{eq:IN2}\end{equation}
 The integrals over $(x_{-},0)$ can be calculated quite similarly.

 It follows that the quantization condition \e{eq:EwrDC}    holds with 
 \begin{align*}
 \varphi_{+}(0, \lambda)&= \lambda_{+}^{1/2+ 1/ \alpha_+} v_+^{-1/ \alpha_+} 
 \alpha_+^{-1}\mathrm{B}(3/2,1/\alpha_+),
  \\
  \varphi_{-}(0, \lambda) & =  \lambda_{-}^{1/2+ 1/ \alpha_-} v_{-}^{-1/ \alpha_-} 
 \alpha_-^{-1}\mathrm{B}(3/2,1/\alpha_-),
\end{align*}
where $\lambda_{\pm}=\lambda-a_{\pm}$.
 In particular, in the case $a_{+}=a_{-}=:a$ the Bohr-Sommerfeld quantization condition reads as
 \begin{align}
( \lambda-a)&^{1/2 + 1/ \alpha_+}  v_+^{-1/ \alpha_+} 
 \alpha_+^{-1}\mathrm{B}(3/2,1/\alpha_+)
 \nonumber \\
   &+ (\lambda-a)^{1/2+ 1/ \alpha_-}  v_-^{-1/ \alpha_-} 
 \alpha_-^{-1}\mathrm{B}(3/2,1/\alpha_-) 
  = \pi\h (n +1/2) + O(\h^2).
\label{eq:INBS}\end{align}
  Plugging  expressions \e{eq:IN1} and \e{eq:IN2} into \e{eq:mainOld} we also find that
 \[
   K_{cl}(\lambda)  = \frac{\lambda_{+}^{1/2+1/\alpha_+} v_+^{-\alpha_+}\alpha_+^{-1}\mathrm{B}(3/2,1/\alpha_+) + \lambda_{-}^{1/2+1/\alpha_-} v_-^{-\alpha_-}\alpha_-^{-1}\mathrm{B}(3/2,1/\alpha_-)}
  {\lambda_{+}^{-1/2+1/\alpha_+} v_{+}^{-\alpha_+}\alpha_+^{-1}\mathrm{B}(1/2,1/\alpha_+) + \lambda_{-}^{-1/2+1/\alpha_-} v_-^{-\alpha_-}\alpha_-^{-1}\mathrm{B}(1/2,1/\alpha_-)}.
  \]
  
Observe that Theorems~\ref{pm} and \ref{norm} can be applied to potential \e{eq:pm1} if $a_{+}=a_{-}$  and $\alpha_{\pm}\geq 2$. If $a_{+}=a_{-}$ but $\alpha_\pm\in [1, 2)$, then we have to use Remark~\ref{discC}; in this case $O(\h^2)$ in \e{eq:INBS} should be replaced by $O(\h^{5/3})$. If $a_{\pm}$ are arbitrary and $\alpha_{\pm}\geq 1$, then the conditions of Theorem~\ref{disc} are satisfied. Moreover, according to Remark~\ref{discR} in the case $\alpha_{\pm}\geq 2$, the estimates of the remainders can be improved. Finally, we note that if $\alpha_{j}<1$, then  $v'' \not\in L^1 ( -\d,  \d)$    so that the semiclassical approximation does not directly work  (even for $a_{+}=a_{-}$) although all formulas above  remain meaningful.
  
\medskip
 
 {\bf 5.4.}
 Let us briefly consider the problem on the half-axis. We now suppose that equation \e{eq:I2} is satisfied for $x\geq 0$, $\psi \in L^2 ({\Bbb R}_{+})$ and
 $\psi (0)=0$. Assumptions~\ref{infi}  and \ref{onewell} should be slightly modified. Namely, we assume that the equation $v(x )=\lambda $ has only one solution $x_{+}=x_{+}(\lambda )$ and $v'(x_{+})>0$ so that $(0, x_{+})$ is a potential well.
 We suppose that the limit of $v(x)$ as $x\to 0$ exists and that the functions $v'(x)$ and 
 $v''(x)$ are bounded in a neighborhood of $x=0$.
 Then the results of Theorem~\ref{Airy} on the solution $u_ + (x)$ of equation \e{eq:I2}  are true  for all $x\geq 0$. In particular, it follows from formula \e{eq:R8-simp} that
  \begin{equation}
 u_+ (0; \lambda,\h)=\pi^{1/2} \h^{1/6} (\lambda  - v(0))^{-1/4}\sin\big( \h^{-1} \int_{0}^{x_{+}(\lambda )} (\lambda  - v(x))^{1/2} dx + \pi/4\big) + O(\h^{7/6}).
 \label{eq:BC}\end{equation}
 Since $\psi (x; \lambda,\h)= c_{+}u_+ (x; \lambda,\h)$, this yields the quantization condition
  \[
\int^{x_{+}(\lambda )}_{0}  (\lambda  - v(x))^{1/2}d x= \pi\h ( n  +3/4) + O(\h^2)
\]
where $n=n (\lambda,\h) $ is an integer. 

Consider now  the boundary condition $\psi' (0)=b \psi(0)$, $b=\bar{b}$. It follows from \e{eq:R8-sim} that
  \begin{align*}
 u_+' (0; \lambda,\h)=-\pi^{1/2} \h^{-5/6} (\lambda  - v(0))^{1/4}\cos\big(   \h^{-1} \int_{0}^{x_{+}(\lambda )} &(\lambda  - v(x))^{1/2} dx + \pi/4\big)
\\
&  + O(\h^{1/6}).
 \end{align*}
Comparing this formula with \e{eq:BC},  we see that the value of $u_+ (0; \lambda,\h)$ is inessential  so that  the quantization condition looks like
  \[
\int^{x_{+}(\lambda )}_{0}  (\lambda  - v(x))^{1/2}d x= \pi\h ( n  +1/4) + O(\h^2).
\]
It does not depend   on $b$.

Other results of Section~4 can also be naturally extended to the problem on the half-axis. 

  \bigskip 
    
I thank D. Robert for a discussion of papers on microlocal analysis.


\begin{thebibliography}{99}
 
 \bibitem {BB} V. M. Babich and  V. S. Buldyrev, {\it Asymptotic methods in   diffraction problems of short-length waves}, Nauka, 1972 (Russian).
  
 
 
  
 
 
   \bibitem {Fed} M. V.  Fedoryuk, {\it  Asymptotic methods for linear ordinary differential equations},  Nauka, 1983 (Russian). 
   
    \bibitem {FedMas} M. V.  Fedoryuk and V. P. Maslov, {\it Semi-classical approximation in quantum mechanics}, Amsterdam, Reidel, 1981.
   

 
 
   \bibitem {ShMaRo} B. Helffer, A. Martinez  et D. Robert, Ergodicit\'e et limite semi-classique,  Comm. Math. Phys., {\bf 109},   313-326, 1987.
   
   
     
   
     \bibitem {ShRo} B. Helffer  et D. Robert, Puits de potentiels g\'en\'eralis\'es,  Ann. Institut H. Poincar\'e, phys. th\'eor., {\bf 41}, No 3, 291-331, 1984.
 
   
   

 
 
\bibitem {LLcl} L. D. Landau and E. M. Lifshitz, {\it Classical mechanics}, Pergamon Press, 1960.


 
 \bibitem {Olver} F. W. J. Olver, {\it Asymptotics and special functions}, Academic
Press, 1974.
 
 
  \bibitem {Simon} B. Simon, Semiclassical analysis of low lying eigenvalues, II, Tunneling, Annals of Math., {\bf 120}, 89-118, 1984.
  
   
   
  

     
    \end{thebibliography}
      \end{document}